\documentclass[11pt]{article}
\usepackage{amsthm}
\usepackage{pxfonts}
\usepackage{yfonts}
\usepackage{dsfont}
\usepackage{graphicx}
\usepackage{relsize}
\usepackage{color}

\def\bel{\begin{equation}\label}
\def\eeq{\end{equation}}
\def\ds{\displaystyle}

\def\mt{\longrightarrow}
\def\v{\vskip 1em}

\def\Rec{{\bf R}}
\def\R{\mathbb R}

\def\C{\mathfrak{B}}

\def\A{{\bf A}}

\def\L{{\bf L}}
\def\U{{\bf U}}

\def\Hat{\widehat}

\def\vol{{\bf vol}}

\def\M{{\bf M}}

\def\Cup{{\bigcup}}

\def\alpha{\alphaup}
\def\beta{\betaup}
\def\gamma{\gammaup}
\def\delta{\deltaup}
\def\theta{\thetaup}
\def\xi{{\xiup}}
\def\eta{{\etaup}}
\def\tau{{\tauup}}
\def\rho{{\rhoup}}
\def\phi{{\phiup}}
\def\psi{{\psiup}}
\def\lambda{{\lambdaup}}
\def\omega{\omegaup}
\def\varphi{{\varphiup}}
\def\gamma{{\gammaup}}

\newtheorem{remark}{Remark}[section]

\begin{document}

\[\hbox{\LARGE{\bf A new proof of maximal theorem on Heisenberg groups}}\]

\[\hbox{Chuhan Sun}\qquad \hbox{and}\qquad \hbox{Zipeng Wang}\]
\begin{abstract}
Given $0\leq\alpha<1$, we define
\[\begin{array}{lr}\ds
\M_\alpha f(u,v,t)~=~\sup_{\Rec\ni(0,0,0)}
\vol \{\Rec\}^{\alpha-1} \iiint_\Rec\left|f [(u,v,t)\odot(\xi,\eta,\tau)^{-1}]\right|d\xi d\eta d\tau
\end{array}\]
where $\Rec\subset\R^{2n+1}$ is a rectangle  parallel to the coordinates. Moreover, $\odot$ denotes the multiplication law on a real Heisenberg group.
\v
The $\L^p$-boundedness of $\M_0$  has been previously proved by M. Christ. 
We show $\M_\alpha\colon\L^p(\R^{2n+1})\mt\L^q(\R^{2n+1})$ for $\alpha={1\over p}-{1\over q},~ 1<p\leq q<\infty$ by applying a geometric covering lemma due to C\'{o}rdoba and Fefferman.
 \end{abstract}

\section{Introduction}
\setcounter{equation}{0}
In this paper, $\odot$ denotes the multiplication law:
\bel{multiplication law}
(u,v,t)\odot(\xi,\eta,\tau)~=~\Big[u+\xi, v+\eta,t+\tau+\mu\big(u\cdot\eta-v\cdot\xi\big)\Big],\qquad\mu\in\R
\eeq
on a real variable representation of Heisenberg groups for every $(u,v,t)\in\R^n\times\R^n\times\R$ and $(\xi,\eta,\tau)^{-1}=(-\xi,-\eta,-\tau)\in\R^n\times\R^n\times\R$.

Let $\Rec\subset\R^{2n+1}$  be a rectangle  parallel to the coordinates.  Given $0\leq\alpha<1$, we define
\bel{M_alpha}
\begin{array}{lr}\ds
\M_\alpha f(u,v,t)~=~\sup_{\Rec\ni(0,0,0)}
\vol \{\Rec\}^{\alpha-1} \iiint_\Rec\left|f [(u,v,t)\odot(\xi,\eta,\tau)^{-1}]\right|d\xi d\eta d\tau.
\end{array}
\eeq
The $\L^p$-boundedness of $\M_0$ has been proved by M. Christ \cite{Michael Christ 2} using  a mixture of techniques developed previously by Ricci and Stein \cite{Ricci-Stein 1}-\cite{Ricci-Stein 2} and Christ \cite{Michael Christ 1} for singular integrals defined on sub-manifolds within a general setting of Nilpotent Lie groups. 

$\diamond$ {\small Throughout, $\C>0$ is regarded as a generic constant depending on its sub-indices.}
\v
{\bf Theorem One}~~{\it Let $\M_\alpha$ defined in  (\ref{M_alpha}) for $0\leq\alpha<1$. We have
\bel{Result One}
\begin{array}{cc}\ds
\left\| \M_\alpha f\right\|_{\L^q(\R^{2n+1})}~\leq~\C_p~\left\| f\right\|_{\L^p(\R^{2n+1})},
\qquad
\alpha~=~{1\over p}-{1\over q},\qquad 1<p\leq q<\infty.
\end{array}
\eeq}

We prove {\bf Theorem One}  by applying a  geometric covering lemma due to C\'{o}rdoba and Fefferman \cite{Cordoba-Fefferman}.

{\bf C\'{o}rdoba-Fefferman covering lemma} ~~
{\it Let $\Rec_j, j=1,2,\ldots$ be a sequence of rectangles in $\R^{2n+1}$ parallel to the coordinates.  There exists a subsequence $\Hat{\Rec}_k, k=1,2,\ldots$ such that
\bel{EST1}
\vol\Bigg\{\Cup_j \Rec_j\Bigg\}~\lesssim~\vol\Bigg\{ \Cup_k \Hat{\Rec}_k\Bigg\}
\eeq
and
\bel{EST2}
\left\|\sum_k \chi_{\Hat{\Rec}_k}\right\|_{\L^p(\R^{2n+1})}~\leq~\C_p~\vol\Bigg\{\Cup_k\Hat{\Rec}_k\Bigg\}^{1\over p},\qquad 1<p<\infty
\eeq
where $\chi$ is an indicator function.}\\
\begin{remark}
This covering lemma is true for any given absolutely continuous measure satisfying the $\A_\infty$-property uniformly on every coordinate subspace. It's proof relies on a beautiful 'slicing' idea which reduces our assertion to a lower dimensional subspace. 

For the Lebesgue measure case, we can prove this result with a more direct approach. 
\end{remark}

 In the next section, we prove  {\bf Theorem One}. In section 3, we give a simpler proof of  the above  covering lemma.

\section{Proof of Theorem One}
\setcounter{equation}{0}
Recall $\M_\alpha$ defined in (\ref{M_alpha}) for $0\leq\alpha<1$. By taking $\xi\mt u-\xi$, $ \eta\mt v-\eta$ and $\tau\mt t-\tau$, $\M_\alpha$ can be equivalently defined as
\bel{M_alpha equi}
\M_\alpha f(u,v,t)~=~
\sup_{\Rec\ni(u,v,t)}\vol\{\Rec\}^{\alpha-1}\iiint_{\Rec}\left|f(\xi,\eta,\tau+\mu(u\cdot\eta-v\cdot\xi))\right|d\xi d\eta d\tau.
\eeq
Let $\alpha={1\over p}-{1\over q}$, $1<p\leq q<\infty$.
Consider
\bel{U_lambda}
\U_\lambda~=~\Bigg\{ (u,v,t)\in\R^{n}\times\R^{n}\times\R\colon \M_\alpha f(u,v,t)>\lambda\Bigg\}.
\eeq
Given any $(u,v,t)\in\U_\lambda$, there is a rectangle $\Rec_j\ni(u,v,t)$ such that
\bel{R_j est alpha}
 \vol\{\Rec_j\}^{\alpha-1}\iiint_{\Rec_j}\left|f(\xi, \eta, \tau+\mu(u\cdot\eta-v\cdot\xi))\right| d\xi d\eta d\tau~>~{1\over 2} \lambda.
\eeq
Let $(u,v,t)$ run through the set $\U_\lambda$. We have
$\ds\U_\lambda\subset\Cup_j~\Rec_j$.
By applying {\bf C\'{o}rdoba-Fefferman covering lemma}, we can select a subsequence $\{\Hat{\Rec}_k\}_{k=1}^\infty$ from the union above and
\bel{union Rk size alpha}
\begin{array}{lr}\ds
\vol\Bigg\{ \U_\lambda\Bigg\}~\leq~\vol\Bigg\{ \Cup_j \Rec_j\Bigg\}
~\lesssim~\vol\Bigg\{ \Cup_k \Hat{\Rec}_k\Bigg\}\qquad \hbox{\small{by (\ref{EST1}) }}.
\end{array}
\eeq
Because $0\leq \alpha<1$,
 we further have
\bel{union size alpha}
\begin{array}{lr}\ds
\vol\Bigg\{ \Cup_k \Hat{\Rec}_k\Bigg\}~\leq~\sum_k\vol \left\{ \Hat{\Rec}_k\right\}
\\\\ \ds
\lesssim~\lambda^{-{1\over 1-\alpha}}\left\{\sum_k \iiint_{\Hat{\Rec}_k}\left|f(\xi, \eta, \tau+\mu(u\cdot\eta-v\cdot\xi))\right| d\xi d\eta d\tau\right\}^{1\over 1-\alpha} \qquad\hbox{\small{by (\ref{R_j est alpha})}}
\\\\ \ds
=~\lambda^{-{1\over 1-\alpha}}\left\{\iiint_{\R^{2n+1}}\left|f(\xi, \eta, \tau+\mu(u\cdot\eta-v\cdot\xi))\sum_k\chi_{\Hat{\Rec}_k}(\xi,\eta,\tau)\right|d\xi d\eta d\tau\right\}^{1\over 1-\alpha}
\\\\ \ds
\leq~\lambda^{-{1\over 1-\alpha}} \left\{ \iiint_{\R^{2n+1}}\left|f(\xi, \eta, \tau+\mu(u\cdot\eta-v\cdot\xi))\right|^p d\xi d\eta d\tau\right\}^{{1\over p} {1\over 1-\alpha}}
\left\|\sum_k\chi_{\Hat{\Rec}_k}\right\|_{\L^{p\over p-1}(\R^{2n+1})}^{1\over 1-\alpha}
\\ \ds~~~~~~~~~~~~~~~~~~~~~~~~~~~~~~~~~~~~~~~~~~~~~~~~~~~~~~~~~~~~~~~~~~~~~~~~~~~~~~~~~~~~~~~~~~~~~~~~~~~
 \hbox{\small{ by H\"older inequality}}
\\ \ds
=~\lambda^{-{1\over 1-\alpha}}  \left\{ \iint_{\R^{2n}}\left\| f(\xi, \eta, \cdot)\right\|_{\L^p(\R)}^p d\xi d\eta \right\}^{{1\over p} {1\over 1-\alpha}} \left\|\sum_k\chi_{\Hat{\Rec}_k}\right\|_{\L^{p\over p-1}(\R^{2n+1})}^{1\over 1-\alpha}
\\\\ \ds
\leq~\big(\C_p\big)^{1\over 1-\alpha}~\lambda^{-{1\over 1-\alpha}}\left\| f\right\|_{\L^p(\R^{2n+1})}^{1\over 1-\alpha} \vol\Bigg\{ \Cup_k \Hat{\Rec}_k\Bigg\}^{{p-1\over p}{1\over 1-\alpha}}
\qquad \hbox{\small{ by (\ref{EST2})}}.
\end{array}
\eeq
By raising both sides of (\ref{union size alpha}) to the $(1-\alpha)$-th power and then taking into account for  $1-\alpha-{p-1\over p}={1\over p}-\left[{1\over p}-{1\over q}\right]={1\over q}$, we find
\bel{union size alpha'}
\vol\Bigg\{ \Cup_k \Hat{\Rec}_k\Bigg\}^{1\over q}~\leq~\C_p ~{1\over \lambda} \left\| f\right\|_{\L^p(\R^{2n+1})}.
\eeq
Let $\U_\lambda$ defined in (\ref{U_lambda}).
From (\ref{union Rk size alpha}) and (\ref{union size alpha'}), we obtain
\bel{weak-type alpha}
\begin{array}{lr}\ds
\vol \Bigg\{ (u,v,t)\in\R^n\times\R^n\times\R\colon \M_\alpha f(u,v,t)>\lambda\Bigg\}^{1\over q}~\leq~\vol\Bigg\{ \Cup_j \Rec_j\Bigg\}^{1\over q}
\\\\ \ds~~~~~~~~~~~~~~~~~~~~~~~~~~~~~~~~~~~~~~~~~~~~~~~~~~~~~~~~~~~~~~~~~~~~~~~~~~~~~~~
~\lesssim~\vol\Bigg\{ \Cup_k \Hat{\Rec}_k\Bigg\}^{1\over q} 
\\\\ \ds~~~~~~~~~~~~~~~~~~~~~~~~~~~~~~~~~~~~~~~~~~~~~~~~~~~~~~~~~~~~~~~~~~~~~~~~~~~~~~~
~\leq~\C_p~{1\over \lambda} \left\|f\right\|_{\L^p(\R^{2n+1})}.
\end{array}
\eeq
By using this weak type $(p,q)$-estimate and applying Marcinkiewicz interpolation theorem, we conclude (\ref{Result One}).

\section{Proof of the covering  lemma}
\setcounter{equation}{0}
Given  $\Rec_j, j=1,2,\ldots$ which are rectangles in $\R^{2n+1}$ parallel to the coordinates, we select $\Hat{\Rec}_k, k=1,2,\ldots$  as follows.

Let $\Hat{\Rec}_1=\Rec_1$. Having chosen $\Hat{\Rec}_1, \Hat{\Rec}_2,\ldots,\Hat{\Rec}_{k-1}$, we pick $\Hat{\Rec}_k$ as the first rectangle $\Rec$ on the list of $\Rec_j$'s after $\Hat{\Rec}_{k-1}$ so that
\bel{R condition}
\vol\left\{~\Rec\cap\bigcup_{\ell=1}^{k-1}\Hat{\Rec}_\ell~\right\}~\leq~\frac{1}{2}\vol\left\{\Rec\right\}.
\eeq
Suppose that $\Rec$ is an unselected rectangle. There is a positive number $N$ such that $\Rec$ is on the list of $\Rec_j$'s after $\Hat{\Rec}_N$ and
\bel{unselect}
\vol\left\{~\Rec\cap\bigcup_{k=1}^N \Hat{\Rec}_k~\right\}~>~\frac{1}{2}\vol\left\{\Rec\right\}.
\eeq
Let $\M$ be the strong maximal operator defined as
\bel{M}
\M f(u,v,t)~=~\sup_{\Rec\ni (u,v,t)} \vol\{\Rec\}^{-1}\iiint_\Rec |f(\xi,\eta,\tau)|d\xi d\eta d\tau.
\eeq
From (\ref{unselect}),  we find
\bel{M>1/2}
\M \chi_{\Cup_k \Hat{\Rec}_k} (u,v,t)~>~{1\over 2},\qquad (u,v,t)\in\Cup_j  \Rec_j.
\eeq
By applying the $\L^p$-boundedness of $\M$, we have
\bel{bound}
\begin{array}{lr}\ds
\vol\Bigg\{\Cup_j\Rec_j \Bigg\}~=~\iiint_{\Cup_j \Rec_j } du dvdt
\\\\ \ds~~~~~~~~~~~~~~~~~~~~~
~\leq~2^2 \iiint_{\Cup_j \Rec_j } \Big[\M \chi_{\Cup_k \Hat{\Rec}_k}\Big]^2 (u,v,t) du dvdt \qquad\hbox{\small{by (\ref{M>1/2})}}
\\\\ \ds~~~~~~~~~~~~~~~~~~~~~
~\lesssim~\iiint_{\R^{2n+1} } \Big[\M \chi_{\Cup_k \Hat{\Rec}_k}\Big]^2 (u,v,t) dudvdt
\\\\ \ds~~~~~~~~~~~~~~~~~~~~~
~\leq~ \iiint_{\R^{2n+1}} \Big[ \chi_{\Cup_k \Hat{\Rec}_k}\Big]^2 (u,v,t) dudvdt
\\\\ \ds~~~~~~~~~~~~~~~~~~~~~
~=~\iiint_{\Cup_k \Hat{\Rec}_k}dudvdt
~=~\vol\Bigg\{\Cup_k\Hat{\Rec}_k \Bigg\}.
\end{array}
\eeq
On the other hand, (\ref{R condition}) suggests
\bel{R condition <}
\vol\left\{~\Hat{\Rec}_k\cap\bigcup_{\ell=1}^{k-1}\Hat{\Rec}_\ell~\right\}~\leq~\frac{1}{2}\vol\left\{\Hat{\Rec}_k\right\},\qquad k=1,2,\ldots
\eeq
which further implies 
\bel{vol E>S}
\vol\Bigg\{ \Hat{\Rec}_k\setminus\Cup_{\ell=1}^{k-1}\Hat{\Rec}_\ell\Bigg\}~>~\frac{1}{2}\vol\left\{\Hat{\Rec}_k\right\},\qquad k=1,2,\ldots.
\eeq
Let $\phi\in\L^{p\over p-1}(\R^{2n+1})$ and $\left\|\phi\right\|_{\L^{p\over p-1}(\R^{2n+1})}=1$. 
We have
\bel{integration}
\begin{array}{lr}\ds
\iiint_{\R^{2n+1}}\phi(\xi,\eta,\tau)\sum_k\chi_{\Hat{\Rec}_k}(\xi,\eta,\tau) d\xi d\eta d\tau
\\\\ \ds
=~\sum_k\iiint_{\Hat{\Rec}_k}\phi(\xi,\eta,\tau) d\xi d\eta d\tau
\\\\ \ds
=~\sum_k \left\{\vol\{\Hat{\Rec}_k\}^{-1}\iiint_{\Hat{\Rec}_k}\phi(\xi,\eta,\tau) d\xi d\eta d\tau\right\} \vol\left\{\Hat{\Rec}_k\right\}
\\\\ \ds
<~2\sum_k \left\{\vol\{\Hat{\Rec}_k\}^{-1}\iiint_{\Hat{\Rec}_k}|\phi(\xi,\eta,\tau)| d\xi d\eta d\tau\right\} \vol\Bigg\{ \Hat{\Rec}_k\setminus\Cup_{\ell=1}^{k-1}\Hat{\Rec}_\ell\Bigg\}\qquad\hbox{\small{by (\ref{vol E>S})}}
\\\\ \ds
\lesssim~\sum_k \int_{\Hat{\Rec}_k\setminus\Cup_{\ell=1}^{k-1}\Hat{\Rec}_\ell} \left\{\vol\{\Hat{\Rec}_k\}^{-1}\iiint_{\Hat{\Rec}_k}|\phi(\xi,\eta,\tau)| d\xi d\eta d\tau\right\} dudvdt
\\\\ \ds
\leq~\sum_k\iiint_{\Hat{\Rec}_k\setminus\Cup_{\ell=1}^{k-1}\Hat{\Rec}_\ell}\M\phi(u,v,t) dudvdt
\\\\ \ds
=~\iiint_{\Cup_k\Hat{\Rec}_k}\M\phi(u,v,t) dudvdt.
\end{array}
\eeq
By using H\"{o}lder inequality and the $\L^p$-boundedness of $M$, we find
\bel{boundedness of M}
\begin{array}{lr}\ds
\iiint_{\Cup_k\Hat{\Rec}_k}\M\phi(u,v,t) dudvdt~\leq~\left\|\M\phi\right\|_{\L^{p\over p-1}(\R^{2n+1})}\vol\Bigg\{\Cup_k\Hat{\Rec}_k\Bigg\}^{1\over p}
\\\\ \ds~~~~~~~~~~~~~~~~~~~~~~~~~~~~~~~~~~~~~~~~~~~~~
~\leq~\C_p~\left\|\phi\right\|_{\L^{p\over p-1}(\R^{2n+1})}~\vol\Bigg\{\Cup_k\Hat{\Rec}_k\Bigg\}^{1\over p}
\\\\ \ds~~~~~~~~~~~~~~~~~~~~~~~~~~~~~~~~~~~~~~~~~~~~~
~=~\C_p~\vol\Bigg\{\Cup_k\Hat{\Rec}_k\Bigg\}^{1\over p}.
\end{array}
\eeq
By substituting (\ref{boundedness of M}) to (\ref{integration}) and taking the supremum of $\phi$, we arrive at
\bel{p norm}
\left\|\sum_k\chi_{\Hat{\Rec}_k}\right\|_{\L^p(\R^{2n+1})}~\leq~\C_p~\vol\Bigg\{\Cup_k\Hat{\Rec}_k\Bigg\}^{\frac{1}{p}},\qquad 1<p<\infty.
\eeq

{\small Westlake University, Hangzhou, 310030, China}\\
{\small email: wangzipeng@westlake.edu.cn}

{\small School of Mathematical Sciences, Zhejiang University, Hangzhou, 310058, China}\\
{\small email: sunchuhan@zju.edu.cn}


\begin{thebibliography}{100}




\bibitem{Cordoba-Fefferman}{\small A.~C\'{o}rdoba and R.~Fefferman, {\it A geometric proof of the strong maximal theorem}, Annals of Mathematics {\bf 102}: 95-100, 1975.}




\bibitem{Michael Christ 1}{\small M.~Christ, {\it Hilbert transforms along curves. I. Nilpotent groups}, Annals of Mathematics {\bf 122}: no.3, 575-596, 1985.}

\bibitem{Michael Christ 2}{\small M.~Christ, {\it The strong maximal function on a nilpotent group}, Transactions of the American Mathematical Society {\bf 331}: no.1, 1-13, 1992.}





\bibitem{Ricci-Stein 1}{\small F.~Ricci and E.~M.~Stein, {\it Oscillatory singular integrals and harmonic analysis on Nilpotent groups}, Proc. Nat. Acad. Sci. U.S.A. {\bf 83}:1-3, 1986.}



\bibitem{Ricci-Stein 2}{\small F.~Ricci and E.~M.~Stein, {\it Harmonic analysis on nilpotent groups and singular integrals. II: Singular kernels
supported on submanifolds}, Journal of Functional Analysis {\bf 78}: 56-84, 1988.}





























\end{thebibliography}
\end{document}